\documentclass[12pt,final,a4paper]{article}%

\usepackage[margin=25.4mm]{geometry}
\usepackage{amsmath,amssymb,graphicx,xcolor}
\usepackage{bm}

\usepackage[utf8]{inputenc}
\usepackage[T1]{fontenc}
\usepackage{bef_alex,tikz}
\def\FG{\bm}
\usepackage[color,notref,notcite,final]{showkeys}
\definecolor{refkey}{gray}{.8}   
\definecolor{labelkey}{rgb}{0.8,0.1,0.4} 

\numberwithin{equation}{section}
\numberwithin{figure}{section}

\newcommand{\ThreeVect}[3]{\!\bma{@{}c@{}} #1\\[-0.1em] #2\\[-0.1em] #3 \ema\!}
\newcommand{\mfC}{\mathfrak C}
\newcommand{\mfU}{\mathfrak U}
\newcommand{\ALPHA}{\rma_{\rml\rmo}}

\begin{document}

\title
{Self-Similar Pattern in \\ Coupled Parabolic Systems as
  \\ Non-Equilibrium Steady States}

\author{
Alexander Mielke%
    \thanks{Weierstra\ss-Institut f\"ur Angewandte Analysis und Stochastik, Berlin,
      Germany  and Humboldt-Universit\"at zu Berlin, Germany.}
\ and 
Stefanie Schindler%
\thanks{Weierstra\ss-Institut f\"ur Angewandte Analysis und Stochastik,
   Berlin, Germany.}
}

\date{January 31, 2023}

\maketitle

\begin{abstract}
  We consider reaction-diffusion systems and other related dissipative systems
  on unbounded domains which would have a Liapunov function (and gradient
  structure) when posed on a finite domain. In this situation, the system may
  reach local equilibrium on a rather fast time scale but the infinite amount
  of mass or energy leads to persistent mass or energy flow for all times. In
  suitably rescaled variables the system converges to a steady state that
  corresponds to asymptotically self-similar behavior in the original system.
\end{abstract}

\section{Introduction}
\label{se:Intro}

Self-similar behavior is a well-studied phenomenon in extended
systems. However, often the view is restricted to simple scalar problems like
the porous medium equation. Moreover, solutions are considered with trivial
behavior at infinity, in particular, in the case of finite mass or energy. 

Here we want to show that a similar behavior occurs in systems of equations but
there we have a richer structure, because pattern may be imposed at
infinity. Rather than looking at systems with traveling pulses or fronts, we
focus on the situation where the local behavior is dominated by a fast trend
towards a unique local equilibrium and the question then arises how the global
solution is evolving through the family of local equilibria. Such phenomena
were studied in \cite{ColEcke90IFES,ColEck92SPSG,vSaHoh92FPSS,EckSch02NSMF} in
the Ginzburg-Landau equation and the Swift-Hohenberg equation. This work is
close to the idea of ``diffusive mixing'' as introduced in \cite{GalMie98DMSS}
for solutions mixing different stable role patterns for $x \to -\infty$ and
$x\to + \infty$.

As the systems under consideration have a ``local gradient structure'', we can
also interpret the self-similar pattern as a \emph{non-equilibrium steady state} and
identify the corresponding fluxes of mass or energy. In particular, we discuss
situations where the scaling leads to a local equilibration of algebraic type
that enforces certain Lagrange multipliers in the diffusive system. In such
cases the Lagrange multipliers can be identified with necessary fluxes that are
needed to understand the mass balances. 

\section{\label{se:IntroPME} The porous medium equation}

As an introduction, we consider the porous medium equation (PME) on the real
line:
\begin{equation}
  \label{eq:PME10}
  u_t = (u^m)_{xx} , \quad t>0,\ x \in \R^1.
\end{equation}
It is well-known that PME has many different self-similar solutions of the form $u(t,x)=
(1{+}t)^{-\alpha} \Phi\big(x/(1{+}t)^\beta\big)$. 

\subsection{\label{su:PME.finite} The finite-mass case} 

The most famous self-similar solution is the Barenblatt solution
\cite{Bare79SSSI} with
\[
\alpha =  \beta= \frac1{m{+}1},\quad W(y)= \max\{0,N-c_my^2\}^{1/(m-1)},
\]
for $m>1$ and $W(y) = N \ee^{-y^2/2}$ for $m=1$,  
where $ N\geq 0 $ can be chosen arbitrary, e.g.\ to achieve the desired total
mass $M=\int_\R W(y) \dd y$. 
This solution can be described as a (non-equilibrium) steady state when
transforming \eqref{eq:PME10} into similarity coordinates. Indeed, setting 
$\tau=\log(1{+}t)$, $y=x/(1{+}t)^{\beta}$, and $w= (1{+}t)^\alpha u$ we find the
equation 
\begin{equation}
  \label{eq:PME12}
  w_\tau = (w^m)_{yy} + \beta yw_y + \alpha w = \big(
      (w^m)_y + \frac1{m{+}1} yw\big)_y.  
\end{equation}
Clearly, $ W(y) = \max\{ N{-}c_my^2,0\}^{1/(m-1)}$ is a steady state, and in
\cite{Vazq07PMEM} there is an extensive study about its global stability. 
\begin{figure}[h]
\begin{minipage}[c]{0.49\textwidth}
\includegraphics[width=\linewidth]{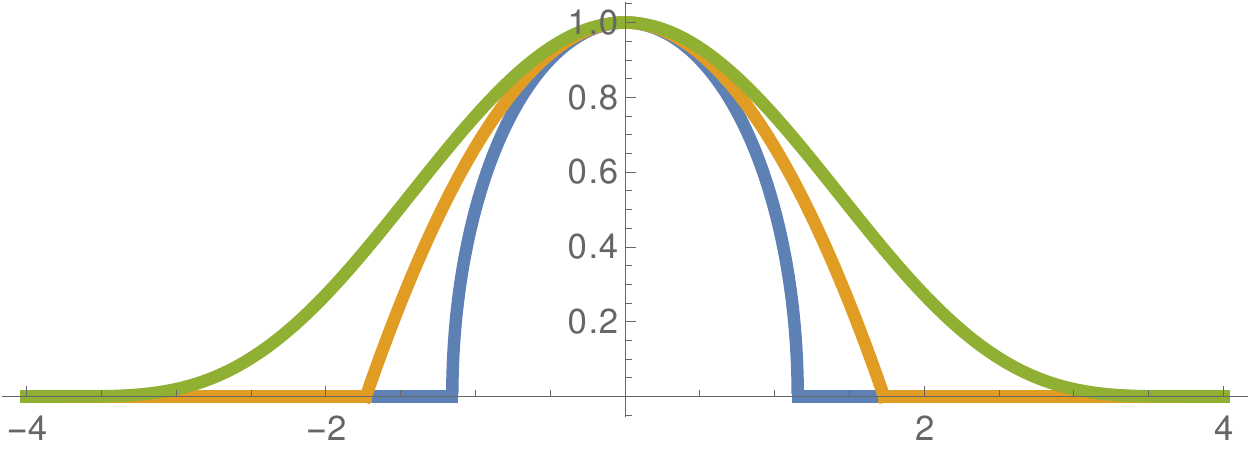}
\end{minipage}
 \begin{minipage}[]{0.49\textwidth}
 \vspace{2.139cm}
  \includegraphics[width=\linewidth]{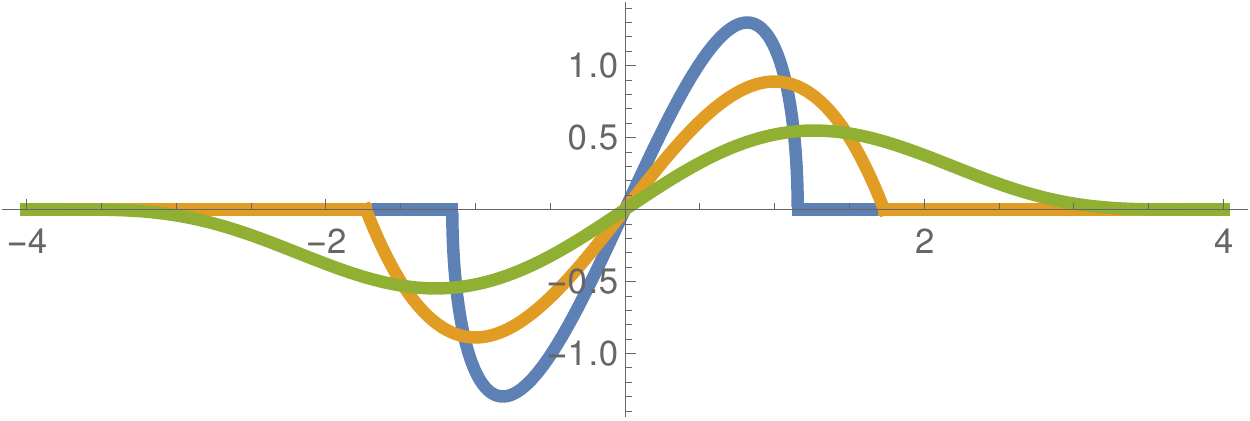}
  \end{minipage}\\
\caption{\label{fig:PME} The left figure shows the Barenblatt profiles $W$ for
  $m=1.25$ (green), $m=2$, and $m=3$ (blue). The right picture shows the
  corresponding self-similar flux pattern $Q$.} 
\end{figure}

To emphasize that $W$ is a non-equilibrium steady state (NESS), we look at the
mass fluxes. In \eqref{eq:PME10} we have the diffusive flux $q(t,x)= -(u^m)_x$
and the total mass $M= \int_\R u(t,x) \dd x $ is preserved for solutions
$u$. Indeed, the flux $q=-(u^m)_x$ takes the form
$q(t,x)= (1{+}t)^{-m\alpha-\beta} Q(x/(1{+}t)^\beta)$ with similarity profile
\[
Q(y) = -\big( W^m\big)_y(y) = -m\,W(y)^{m-1}  W'(y) . 
\]

\subsection{Diffusive mixing and infiltration}

We may also consider PME with boundary conditions $u(t,\pm\infty)= U_\pm$ with
different concentrations $U_-$ and $U_+$. Again a self-similar profile develops
but now the scaling is different as $u$ cannot be scaled by a prefactor,
because of the boundary conditions. The boundary conditions provide reservoirs
with an infinite amount of mass at $x=+\infty$ and $x = - \infty$. The
diffusive mixing describes how the mass is flowing from one reservoir to the
other. 

\begin{figure}[h] 
\includegraphics[width=0.49\linewidth]{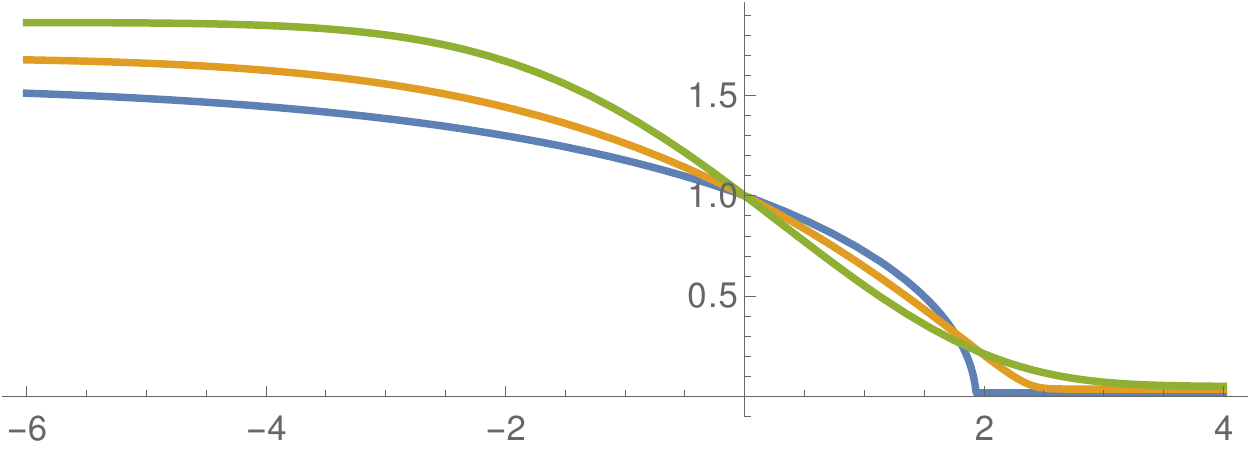}
\includegraphics[width=0.49\linewidth]{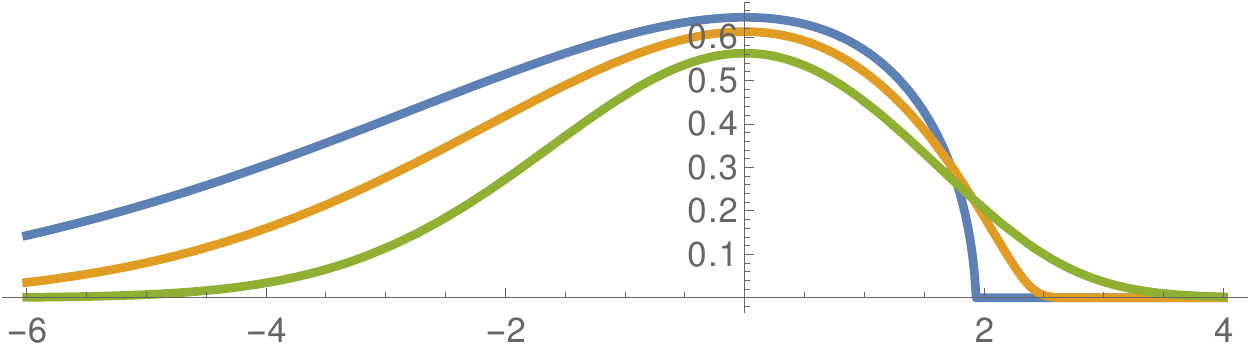}
\caption{\label{fig:PME.mixing} The left figure shows the self-similar infiltration
  profiles with $U_+=0$ and $W(0)=1$ for $m=1.25$ (green), $m=2$, and $m=3$
  (blue). The right picture shows the corresponding flux $Q$.} 
\end{figure}

In particular, we have to choose $\alpha=0$ and are then forced
to take $\beta=1/2$, which is the parabolic scaling. The corresponding equation
in the parabolic similarity coordinates reads
\begin{equation}
  \label{eq:PME.30}
  w_\tau = (w^m)_{yy} + \frac y2\, w_y, \quad w(\tau,\pm\infty)=U_\pm.
\end{equation}
Of course steady states $W$ are again exact self-similar solutions to
\eqref{eq:PME.30}.  
The existence and uniqueness of stationary profiles $W$ with
$W(\pm\infty)=U_\pm$ are studied in
\cite{MieSch21?ESPS}. The profiles are monotone and converge to their limits
$U_\pm$ faster than exponential. For $U_- > U_+$ the flux is nonnegative and
has its maximum at $y=0$, see Figure \ref{fig:PME.mixing}. The diffusive flux
$q=-(u^m)_x$ scales differently from before, but the self-similar profile $Q$ has
the same expression as before:
\[
\begin{aligned}
&q(t,x)=- (1{+}t)^{-1/2} Q\big(x/(1{+}t)^{1/2}\big)  \ \text{ with }
\\
&Q(y) = -\big( W^m\big)'(y) = - m W(y)^{m-1} W'(y).
\end{aligned}
\]

The case $U_+=0$ is called the case of filtration, where mass is flowing into
the area $x \gg 1$ where initially the concentration is $0$. For all $m>1$ the
front propagates like $t^{1/2}$ and the infiltrated mass given by
$M_\geq (t) = \int_0^\infty u(t,x)\dd x $ satisfies $\dot M_\geq (t) =
q(t,0)=(1{+}t)^{-1/2} Q(0)$, i.e.\ we have $M_\geq (t) = M_\geq (0)(1{+}t)^{1/2}$.

\section{\label{se:CGPS} A model motivated by turbulence}

Kolmogorov's two-equation model \cite{Kolm42ETMI,Spal91KTEM,BulMal19LDAK}
considered on all of $\R^d$ has a rich scaling structure and hence allows for
self-similar solutions, see Sec.\,3 in \cite{MieNau22EGTW}. In
\cite{Miel21?TCDP} a simplified model is studied, where $\wt v(t,x)$ is a
scalar shear velocity and $\wt k(t,x)$ is the mean turbulent kinetic energy:
\begin{align*}
\wt v_t &= \DIV\big( \eta \,\wt k^\beta \nabla v\big), 
\\
\wt k_t &= \DIV\big(\kappa\, \wt k^\beta \nabla \wt k\big) 
   + \eta \wt k^\alpha |\nabla \wt v|^2,
\end{align*}
where $\eta,\kappa, \beta$ are positive parameters. Note that the system
contains the PME with $m=\beta+1$, if we look at the case $\wt v\equiv 0$. 

The system has the total linear momentum $\calP(\wt v)=\int_{\R^d} \wt v(x) \dd x$ and
the total energy $\calE(\wt v,k) =\int_{\R^d} \big( \frac12 \wt v^2 + \wt k\big) \dd x$ as
conserved quantities. The kinetic energy that is dissipated via shear viscosity
(depending on $\wt k$) is fully fed into the turbulent kinetic energy, which
leads to the energy conservation. 

In fact, the system can be written as a gradient-flow equation with respect to
entropy $\calS(\wt k)= \int_{\R^d} \wt k^\theta \dd x $ for any
$\theta\in (0,1)$, hence it is expected that $\wt k$ has to become
constant. For bounded domains with no-flux boundary conditions, it can be shown
that solutions converge exponentially to the unique equilibrium state with
constant $\wt v$ and $\wt k$ such that the conserved quantities match, see
Sec.\,2 in \cite{Miel21?TCDP}.

On the unbounded domain $\R^d$, solutions with finite momentum $\calP$ and
finite energy $\calE$ are expected to disperse and converge uniformly to $0$.
In Sec.\,6 of \cite{Miel21?TCDP} it is argued that $\wt v$ and $\wt k$ behave
asymptotically self-similar, but with different exponents because of
$\frac12 \wt v^2$ in the energy and $\wt v$ in the momentum. The conjecture is
that $\wt k$ develops, for large $t$, a self-similar pattern of Barenblatt
type, namely 
\begin{equation}
  \label{eq:Barenblatt}
  K(y) = \max\{0, N-c_{\beta+1}|y|^2\}^{1/\beta}
\end{equation} 
with total mass $E=\calE(v(0),k(0))$,
which means that all macroscopic kinetic energy is converted into turbulent
kinetic energy. Moreover, $\wt v$ develops, for large $t$, a self-similar
pattern that is a Barenblatt solution
raised to the power $\kappa/\eta$, i.e.\ 
\[
V(y)= \wt v\, \big(K(y)\big)^{\kappa/\eta}
\]
with $\wt v $ such that $\calP(V)= \calP(v(0))$.  

More precisely, with $\gamma=1/(2{+}\beta d)$ we rescale the variables via
\[
\tau=\log(1{+}t),  \ y= x/(1{+}t)^\gamma, \ k=(1{+}t)^{\gamma d} \wt k, \
v=(1{+}t)^{\gamma d} \wt v
\]
and obtain a non-autonomous coupled system 
\begin{align*} 
v_\tau& = \DIV\big( \gamma v y + \eta\, k^\beta \nabla v\big), \\
k_\tau&= \DIV\big( \gamma k y + \kappa\, k^\beta \nabla k\big) 
  + \ee^{-\gamma d \tau} \, \eta k^\beta |\nabla v|^2, 
\end{align*}
where now $\DIV$ and $\nabla$ are taken with respect to $y$.  Thus, we see that
the equation for $k$ behaves, for $\tau \gg 1$, like the PME and has the
Barenblatt profiles from \eqref{eq:Barenblatt} as asymptotic steady
states. Inserting such a Barenblatt solution $K$, one can show that the linear
equation for $V$ has the steady states $\wt v K^{\kappa/\eta}$. For $\eta \gg
\kappa$, this means that $V$ will have large gradients near the boundary
of the support of $K$.\medskip

As for the PME equation there are also solutions with infinite momentum or
energy, because there are nontrivial limits at $x \to \pm \infty$.  In one such
case it is possible to write down an exact self-similar solution, namely for
$\alpha=\beta=1$ and $\kappa=\eta$, where we set $\eta=1$ for simplicity.  With
$y=x/\sqrt{1{+}t}$ and arbitrary $A>0$ it can be checked that $ v(t,x)=V(y)$
and $k(t,x)= K(y)$ are explicit solutions, if we choose
\[
\big(V(y),K(y)\big)=\left\{ \ba{@{}cl} \big(A/\sqrt{2},0\big) &\text{for }
  y\geq A, \\ \big(y/\sqrt{2}, (A^2{-}y^2)/4\big) &\text{for }|y|\leq A, \\ 
\big(-A/\sqrt{2},0\big)& \text{for } y\leq -A . \ea\right.
\]
For this solution, the energy density $e(t,x)=\frac12v(t,x)^2 + k(t,x)$ is
indeed equal to the constant $A^2/4$, which means that the solutions have
infinite total energy. Nevertheless there are nontrivial fluxes, namely for the
linear momentum and the turbulent kinetic energy:
\begin{align*}
&Q^\mafo{lin.mom}(y) = -K(y)V'(y)= -K(y)/\sqrt{2}, 
\\
&Q^\mafo{tur.kin}(y) = -K(y)K'(y), \text{ and}
\\
&S^\mafo{tur.kin}(y) = K(y) V'(y)^2 =  K (y)/2\geq 0,
\end{align*}
where the last term is the source of turbulent kinetic energy stemming from the
dissipation in the momentum equation.

\section{\label{se:RDS} Diffusive mixing in reaction-diffusion systems}

Here we consider systems of equations which describe the concentrations $c_j$
of species $X_j$ that diffuse with a diffusion constant $d_j$ and that undergo
reactions according to the mass-action law. Our main assumption is that there
is a continuous family of equilibria to the reaction equation
$\dot \bfc = \bfR(\bfc)$, where
$\bfc=(c_1,...,c_{i_*}) \in \mfC:={[0,\infty[}^{i_*} $ is the vector of
concentration.  We consider the reaction-diffusion system (RDS) on the real
line $\R^1$ and impose boundaries conditions at infinity which represent
reservoirs of infinite mass. The prescribed limit states $C_-$ and $C_+$ at
$x=\pm \infty$ are assumed to be in equilibrium, i.e.\ $\bfR(C_\pm)=0$. Thus,
our RDS takes the form
\begin{equation}
  \label{eq:unscRDS}
    \wt\bfc_t= \bfD\, \wt\bfc_{xx} + \bfR(\wt\bfc),  \quad \wt \bfc(t,\pm\infty)=C_\pm,
\end{equation}
where $\bfD$ is the diagonal matrix $\mafo{diag}(d_1,...,d_{i_*})$. 

To study self-similar behavior, we use the parabolic scaling variables
$\tau = \log(1{+}t)$ and $y= x/\sqrt{1{+}t}$ again and find for
$\bfc(\tau,y)=\wt\bfc(t,x)$ the scaled equations
\begin{equation}
  \label{eq:scRDS}
  \bfc_\tau = \bfD \, \bfc_{yy} + \frac y2 \, \bfc + \ee^\tau \,\bfR(\bfc),\quad 
  \bfc(\tau,\pm\infty)=C_\pm, 
\end{equation}
where the prefactor $\ee^\tau$ appears because the reactions do not scale in a
similar way as the derivatives $\pl_\tau$ and $\pl_y^2$.
 
Hence, for large $\tau$ the reaction becomes stronger and stronger and will lead
to a local equilibration of the reactions. Of course, this is only an effect of
the scaling, but it says that on long time scales we first see that the
reactions act on their natural time scale while the diffusive mixing may take
much longer and will actually never stop because of the boundary conditions at
$\pm \infty$.

\subsection{\label{su:LargeTimeLim} The diffusive large-time limit and
  reduced systems}

We can formally go to the limit $\tau\gg 1$ in \eqref{eq:scRDS} as
follows. Clearly, the reaction has to become equilibrated, i.e.\
$\bfR(\bfc(\tau,y))=\bm0$ for all $\tau$ and $y$. However, the product $\ee^\tau
\bfR(\bfc)$ should then be treated as a limit of the type ``$\infty \cdot
\bm0$'', taking the value $\bflambda(\tau,y)\in \R^{i_*}$. This term can be
understood as a rescaled version of a small reaction flux, because the reaction
$\bfc$ will only be equilibrated up to order $\ee^{-\tau}$ such that
$\bfR(\bfc)$ may still contain a term $\ee^{-\tau} \bflambda $. 

The models resulting from \eqref{eq:unscRDS} and \eqref{eq:scRDS} are the
following constrained RDS 
\begin{subequations}
  \label{eq:unscscConstrRDS}
\begin{align}
  \label{eq:unscConstrRDS}
  \wt\bfc_\tau &= \bfD\,\wt\bfc_{xx} +\wt\bflambda, \quad \bfR(\wt\bfc)=0, \quad
  \wt\bfc(t,\pm\infty) = C_\pm,
\\
  \label{eq:scConstrRDS}
  \bfc_\tau &= \bfD\,\bfc_{yy} +\frac y2\,\bfc_y +\bflambda, \quad \bfR(\bfc)=0, \quad
  \bfc(\tau,\pm\infty) = C_\pm. 
\end{align}
\end{subequations}
The important point is that $\bflambda$ is restricted to lie in the linear
subspace
$\mafo{span}\bigset{\rmD\bfR(\bfc)\bfv}{ \bfR(\bfc)=0,\ \bfv\in \R^{i_*}}$,
such that $\bflambda$ plays the role of a Lagrange multiplier to the constraint
$\bfR(\bfc)=\bm 0$.

We restrict to the case of mass-action kinetics, where the reaction is in
detailed balance. Then, there exists a surjective linear stoichiometric mapping
$\bfQ:\R^{i_*} \to \R^{j_*} $ giving the conserved molecular masses and a
nonlinear map $\Psi : \;\mfU:=\bfQ\mfC \to \mfC$ such that
\begin{align*}
&\bfQ\bfR(\bfc)=0  \text{ for all } \bfc \in \mfC, \quad
\bfQ\Psi(\bfu)=\bfu  \text{ for all } \bfu \in \mfU, 
\\
&\text{and }\bigset{\bfc \in \mfC }{ \bfR(\bfc)=0} = \bigset{\psi(\bfu) \in \mfC}{
  \bfu \in \mfU}, 
\end{align*}
i.e.\ $\Psi$ parametrizes the set of equilibria of $\bfR$. We refer to 
\cite{MiPeSt21EDPC,MieSch22?CSSP} for more details. 

Setting $\wt\bfu(t,x)=\bfQ \wt\bfc(t,x)$, $\bfu(\tau,y)=\bfQ \bfc(\tau,y)$, and
$U_\pm = \bfQ C_\pm$, the constrained RDS \eqref{eq:unscscConstrRDS} in unscaled and
scaled form reduce to the simple diffusion systems
\begin{subequations}
  \label{eq:unscscDiff.bfu}
\begin{align}
  \label{eq:unscDiff.bfu}
\wt\bfu_t &= \big( \bfA(\wt\bfu)\big)_{xx}  , \quad  \bfu(t,\pm\infty) = U_\pm
\quad \text{and}
\\
  \label{eq:scDiff.bfu}
\bfu_\tau &= \big( \bfA(\bfu)\big)_{yy} + \frac y2\, \bfu_y , \quad  \bfu(\tau,
\pm\infty) = U_\pm 
\end{align}
\end{subequations}
with $\bfA(\bfu) = \bfQ\, \bfD\, \Psi(\bfu)$. 
Note that $\bfQ \bflambda\equiv 0$ by construction. 

\subsection{\label{su:VectValProfileEqn} Vector-valued profile equations}

Under the assumption that $\bfA:\R^{j_*} \to \R^{j_*}$ is (strongly) monotone, the
vector-valued profile equation
\begin{equation}
  \label{eq:Profile.bfU}
  \bm0= \big( \bfA(\bfU)\big)_{yy} + \frac y2\,\bfU_y , \quad \bfU(\pm\infty)=U_\pm
\end{equation}
has a unique similarity profile $\bfU:\R\to \R^{j_*}$. We refer to
\cite{GalMie98DMSS} for the scalar-valued case and to \cite{MieSch21?ESPS} for
the more general vector-valued case. 

A profile $\bfU$ solving \eqref{eq:Profile.bfU} is a classical steady state
solution for the scaled diffusion system \eqref{eq:scDiff.bfu}. Hence, setting 
$\wt u(t,x)= \bfU\big(x/\sqrt{1{+}t}\big)$ provides an exact self-similar
solution to \eqref{eq:unscDiff.bfu}.

Clearly, defining $\bfC(y) = \Psi(\bfU(y))$ we obtain a solution $\bfC:\R\to
\R^{i_*}$ for the constrained profile equation 
\begin{equation}
  \label{eq:Profile.bfu}
\begin{aligned}
  &\bm0= \bfD \bfC_{yy} + \frac y2\,\bfC_y + \bflambda , \quad \bfR(\bfC)=\bm0, 
\\
&\bfQ\bflambda =0, \quad \bfC(\pm\infty)=\Psi(U_\pm). 
\end{aligned}
\end{equation}
The similarity profile $\bfC$ is a steady state for the scaled constrained RDS
\eqref{eq:scConstrRDS}, and $\wt\bfc(t,x)=\bfC\big(x/\sqrt{1{+}t}\big) $ is an
exact self-similar solution for \eqref{eq:unscConstrRDS}.\medskip

In the following three subsections, we consider a few special cases, where we
highlight the role of the reaction flux(es) $\bflambda$ in particular.

\subsection{One reaction for two species}
\label{su:OneReact}

In \cite{GalSli22DREE, MieSch22?CSSP} the following system of two equations is
studied in detail:
\[
\binom{\dot c_1}{\dot c_2} = \binom{d_1\, \pl_x^2c_1}{d_2 \,\pl_x^2 c_2} 
 +  \kappa \big( c_2^\beta {-}  c_1^\gamma\big) \binom{\gamma}{-\beta } \qquad \text{for } t>0 \text{ and } x \in \R. 
\]
The two concentrations $c_1, c_2\geq 0$ for the species $X_1,X_2$ 
diffusive with diffusion constants $d_j$ and
undergo the reversible mass-action reaction pair $\gamma X_1 \rightleftharpoons
\beta X_2$.   

The scaled and constraint system \eqref{eq:scConstrRDS} 
takes the form
\[
\pl_\tau\! \binom{c_1}{c_2} = \binom{d_1\, \pl_y^2c_1}{d_2\, \pl_y^2 c_2} + \frac y2
\pl_y\!\binom{c_1}{c_2} + \Lambda \binom{\gamma}{\!\!-\beta\!\!}, \ \ \Lambda \in \R,
\ \  c_1^\gamma = c_2^\beta .
\]
Here $\bflambda = \Lambda(\gamma,-\beta)^\top \in \R^2$ contains
only one scalar reaction
flux $\Lambda \in \R$, because there is only one reaction pair.
 
The set of equilibria for $\bfR$ is the one-parameter family
\[
\bigset{\bfc\in \mfC}{ \bfR(\bfc)=0} = \bigset{(A^\beta,A^\gamma)}{ A\geq 0}.
\]
The linear stoichiometric mapping is
$\bfQ=\big( \beta\ \ \gamma\big) \in \R^{1\ti 2}$ defining
$u=\bfQ\bfc=\beta c_1 {+} \gamma c_2\geq 0$, and $\Psi:{[0,\infty[}=\mfU\to \mfC$
is defined via
\[
\bfc=\Psi(u)=\binom{\psi_1(u)}{\psi_2(u)}    \quad
\Longleftrightarrow \quad \left\{ \ba{c} u=\bfQ\bfc=\beta c_1{+}\gamma c_2\\
 \text{and } c_1^\gamma  = c_2^\beta \ea \right.
\]
The case $ \gamma = \beta $ leads to the simple relation $\Psi(u) = 
\frac1{\beta{+}\gamma} \binom uu$.  
If $ \beta \neq \gamma $ we may assume $\beta < \gamma$ without loss of
generality, see \eqref{eq:ExamplePsi}  for a nontrivial example.

For $A_\Psi(u):= \bfQ\bfD \Psi(u)=\binom{\beta d_1}{\gamma d_2} \vdot \Psi(u)$
one obtains $0<\psi'_1 (u) \leq \psi'_1(0)=1/\beta$ and
$0<\psi'_2(u) \leq \psi'_2(\infty) = 1/\gamma$. This yields
\[
D_*=\min\{d_1,d_2\} \leq A'_\Psi(u) \leq D^*=\max\{d_1,d_2\}
\]       
as well as $A'_\Psi(u)\to d_1 $ for $u\to 0^+$ and 
$A'_\Psi(u) \to d_2$ for $u\to \infty $. 

Thus, the existence theory for similarity profiles in
\cite{GalMie98DMSS,MieSch21?ESPS} provides a unique and smooth solution $U$ of the 
profile equation 
\[
\big(A_\Psi(U)\big)'' + \frac y2 \, U'=0 \ \text{ on } \R, \quad U(\pm\infty)=U_\pm.
\]
Assuming $U_-<U_+$, this solution is strictly increasing and converges to its two
limits like the error function. In addition to $U_- < U(y) < U_+$ the estimate
\begin{equation}
  \label{eq:U'Estim}
  0 < U'(y) \leq \ee^{-y^2/(4D^*)} \sqrt{\tfrac{\ds D^*}{\ds 8D_*^2}}\, \big( U_+ -
U_-\big) \quad \text{for all }y \in \R
\end{equation}
holds, even in the
case $U_-=0$, where asymptotically the concentrations vanish, viz.\ 
$C_-=\Psi(U_-)=\binom00$, because the
effective diffusion is still bounded from below by $D_*>0$.    

Such a profile $U:\R\to [U_-,U_+]$ for the reduced equation leads to a smooth
concentration profile $\bfC:\R\to \mfC\subset \R^2$ given by
$\bfC(y)=\Psi(U(y))$ and satisfying the profile equation
\begin{align*}
&0= \bma{@{}c@{}c@{}}d_1&0\\ 0&d_2\ema \bfC''+ \frac y2 \,\bfC' + \Lambda
\binom\gamma{-\beta}, \quad C_1^\gamma=C_2^\beta, \\
&
 \bfC(y)\to \Psi(U_\pm) \text{ for }y\to \pm\infty.
\end{align*}
Hence, the reaction flux $\Lambda$ can be written as 
\[
\Lambda(y):= -\frac1{\gamma}\big(d_1 C''_1(y) + \frac y2 C'_1(y)\big) =
\frac1{\beta }\big(d_2 C''_2(y) + \frac y2 C'_2(y)\big). 
\]

In general, $\Lambda$ will be nontrivial, this can already be seen in the
simple case $\beta=\gamma$, which implies $C_1\equiv C_2$,
$\psi_j(u)=u/(\beta{+}\gamma)$, and hence $A_\Psi(u)=\frac{d_1{+}d_2}2
u$. Denoting by $\bbE:\R\to {]0,1[}$ the unique solution of $\bbE'' + y
\bbE'=0$, $\bbE(-\infty)=0$, and $\bbE(\infty)=1$ and recalling $U_\pm =
2\gamma C_\pm$, we obtain the unique profiles 
\begin{align*}
U(y) &= 2\gamma C_- + 2\gamma(C_+{-}C_-) \bbE\big(y/(d_1{+}d_2) \big),
\\
C_1(y)&=C_2(y)= \frac1{2\gamma} \,U(y).
\end{align*}
This provide the explicit formula (for $\beta=\gamma$ only), namely
\[
\Lambda(y) = \frac{2(d_1{-}d_2)}{(d_1{+}d_2)^2} \,\big( C_+-C_-\big)
\:\bbE''\big( y/(d_1{+}d_2)\big),
\]
i.e.\ only for $d_1=d_2$ we have $\Lambda \equiv 0$. 

\begin{figure}
\centering \begin{tikzpicture}
\node at (0,0){\includegraphics[width=0.7\linewidth]{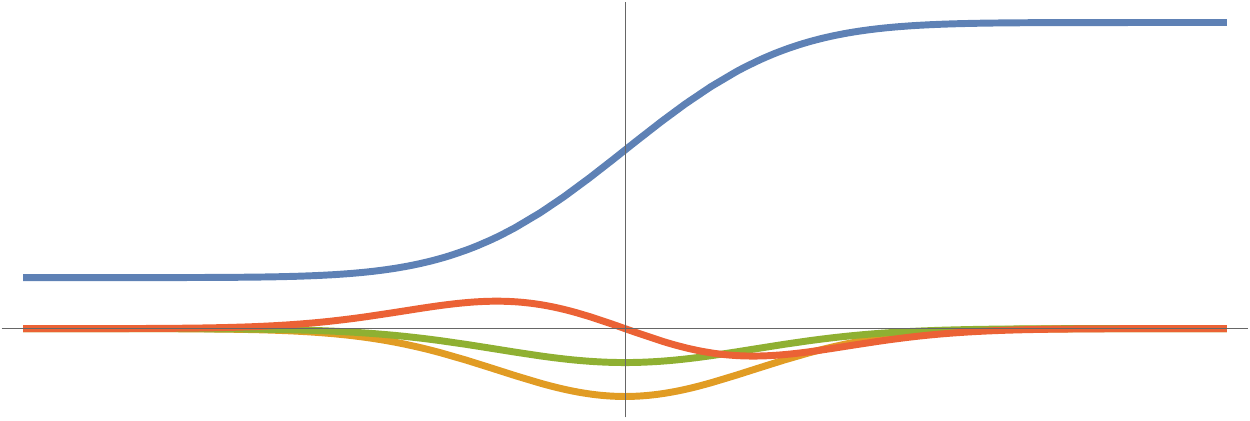}};
\node at (5.5,-1.3){$y$};
\node[color=blue!60!black] at (3.7,1.3){$C_1=C_2$};
\node[color=brown!50!orange] at (-1.7,-1.8){$Q_1^\mafo{diff}$};
\node[color=green!30!brown] at (1.6,-0.7){$Q_2^\mafo{diff}$};
\node[color=red!60!orange] at (-0.55,-0.5){$\bm\Lambda^\mafo{react}$};
\end{tikzpicture}
\caption{\label{fig:TwoOneRDS11} For the case $\beta=\gamma=1$ the similarity
  profile $C_1=C_2$ is shown together with the diffusive fluxes
  $Q_j^\mafo{react}$ and the reaction flux $\Lambda^\mafo{react}$.} 
\end{figure}
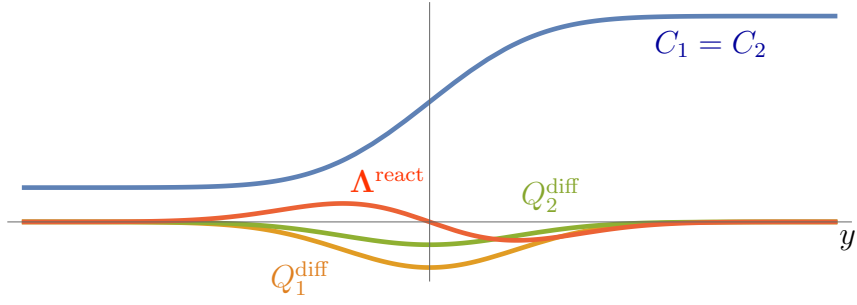

In Figure \ref{fig:TwoOneRDS11} we display, for $C_-=0.2<
C_+=1.2$, $\beta=\gamma=1$, and $d_1=1>d_2=0.5$,  the profile $C_1=C_2$, 
the associated diffusion fluxes
$Q^\text{diff}_j= -d_j C'_j(y)$ for $j=1,2$, and the reaction flux
$\Lambda^\mafo{react}$. Because of the $d_1>d_2$ the 
diffusive fluxes satisfy $|Q^\text{diff}_1|> |Q^\text{diff}_2|$, so one would
expect the profile $C_1$ to be flatter than $C_2$. However,  $C_1=C_2$ is
realized by the reaction $X_1 \rightleftharpoons X_2$, which pushes missing or
excessive mass from $X_1$ into $X_2$.  

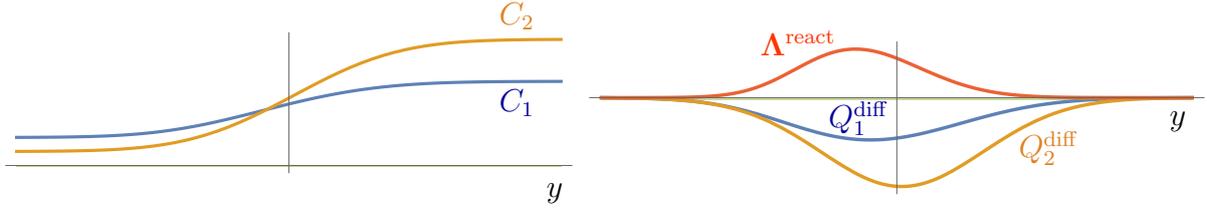
\begin{figure}[h]
\centering 
\begin{tikzpicture}
\node at (-4,0){\includegraphics[width=0.47\linewidth]{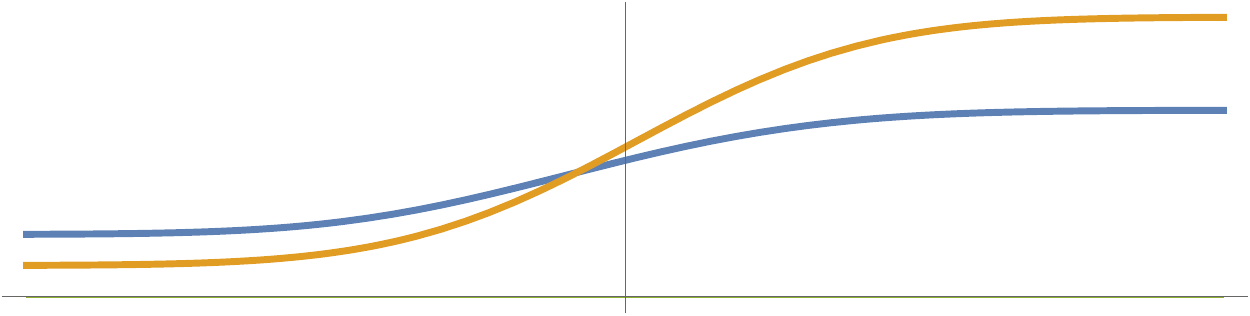}};
\node at (4,-0.2){\includegraphics[width=0.51\linewidth]{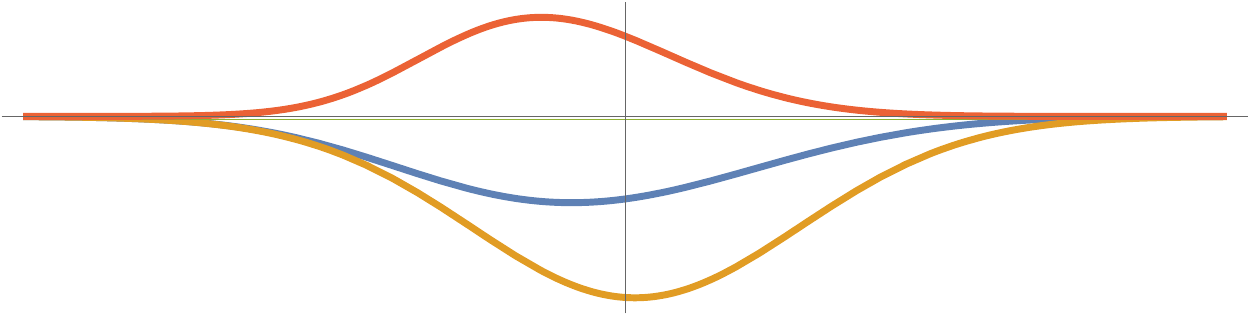}};
\node at (-0.5,-1.2){$y$};
\node at (7.7,-0.25){$y$};
\node[color=blue!60!black] at (-1,0){$C_1$};
\node[color=brown!50!orange] at (-1,1.15){$C_2$};
\node[color=blue!60!black] at (3.5,-0.2){$Q_1^\mafo{diff}$};
\node[color=brown!50!orange] at (6,-0.6){$Q_2^\mafo{diff}$};
\node[color=red!60!orange] at (2.7,0.8){$\bm\Lambda^\mafo{react}$};
\end{tikzpicture}


\caption{\label{fig:RDSTwoOne12} For the case $\beta=1<\gamma=2$ and
  $d_1=d_2=1$ the similarity
  profile $C_1$ and $C_2$ are shown (left picture) and the associated
  diffusion fluxes and reaction flux (right picture).} 
\end{figure}

We also consider the case $\beta=1$ and $\gamma=2$ which corresponds to the
nonlinear reaction pair $X_1  \rightleftharpoons 2X_2$. Now, the profiles are
no longer identical and there is no symmetry $y\leftrightarrow -y$. We find
\
\begin{equation}
  \label{eq:ExamplePsi}
  \Psi(u)=\binom{ \frac14\big(\sqrt{1{+}8u} -1\big)}
                {\frac18\big(1+4u- \sqrt{1{+}8u} \big) },
\end{equation}
and can calculate all fluxes for $C_-=\Psi(1)=\binom{1/2}{1/4}$ and
  $C_-=\Psi(6)=\binom{3/2}{9/4}$, now choosing $d_1=1$ and $d_2=1$ which gives
  $A_\Psi(u)=u$ and makes the calculation simple. We refer to Figure
  \ref{fig:RDSTwoOne12} for the corresponding profiles and diffusion and
  reaction fluxes.  

\subsection{One reaction for three species}
\label{su:OneReact3Spec}

For the typical binary reaction $X_3 \rightleftharpoons  X_1 {+} X_2$ we 
obtain the scaled constrained {\setlength{\arraycolsep}{0.12em} RDS 
\[
\pl_\tau \bfc 
 = \bfD \pl_{y}^2 \bfc   
   + \frac y2  \pl_y \bfc 
  +\Lambda \ThreeVect11{-1}, \quad c_1c_2=c_3 
\]
with} $\bfD=\mafo{diag}(d_1,d_2,d_3)$. The profile equation reads
\begin{equation}
  \label{eq:ProfEqn3Spec}
\begin{aligned}
0 &= \bfD \bfC'' + \frac y2  \bfC' + \Lambda \ThreeVect11{-1} , \quad 
    C_1C_2=C_3 
\quad 
\text{and} \  \     \bfC(\pm\infty) = \Psi(\bfU_\pm). 
\end{aligned}
\end{equation}

The set of equilibria for $\bfR$ is a two-parameter family:
\[
\bigset{\bfc\in \mfC}{ \bfR(\bfc)=0} = \bigset{(A,B,AB)}{ A, B\geq 0}.
\]
We can choose the stoichiometric matrix 
\[
\bfQ = \bma{ccc} 1&0&1\\ 0&1&1\ema \in \R^{2\ti 3} 
\]
and obtain $\bfu=\binom{u_1}{u_2} = \bfQ\bfc \in \mfU:={[0,\infty[}^2$. The reduction
function $\Psi:\mfU\to \mfC$ can be calculated explicitly in the form
\begin{align*}
\Psi(u_1,u_2)&=\frac12
  \ThreeVect{ u_1{-}u_2{-}1 + s(\bfu) }
            { u_2{-}u_1{-}1 + s(\bfu) }
            { u_1{+}u_2{+}1 - s(\bfu) } \quad \text{with}
 s(\bfu)=\sqrt{(1{+}u_1{+}u_2)^2-4u_1u_2}.
\end{align*}
To extend $s$ to a function $s:\R^2\to \R$ we simply set $s(u_1,u_2)=1+u_1+u_2$
whenever $u_1\leq 0$ or $u_2\leq 0$ and observe that $s$ is globally Lipschitz
continuous. Moreover, $s_j(\bfu)=\pl_{u_j} s(\bfu)$ satisfies $s_1(\bfu)\leq
1$, $s_2(\bfu)\leq 1$ and $s_1(\bfu) +s_2(\bfu)\geq 0$ for all $\bfu\in \R^2$. 

From this we can calculate the function $\bfA(\bfu)=\bfQ\,\bfD\,\Psi(\bfu)$:
\[
\bfA(\bfu)= \frac12 \binom{(d_1{+}d_3) u_1 + (d_3{-}d_1)(1{+}u_2{-}s(\bfu)) }
  { (d_2{+}d_3) u_2 + (d_3{-}d_2)(1{+}u_1{-}s(\bfu)) } .
\]
To show monotonicity of $\bfA:\R^2\to \R^2$ we observe that 
for general $\rmC^1$ functions $\bfA$ we have the equivalence
\begin{align*}
&\forall\, \bfu,\wt\bfu:\ \langle\bfA(\bfu){-}\bfA(\wt\bfu), \bfu{-}\wt\bfu \rangle \geq
\ALPHA |\bfu{-}\wt\bfu |^2 
\\
&\Longleftrightarrow \quad 
\forall\, \bfu:\ \frac12\big(\rmD\bfA(\bfu){+}\rmD\bfA(\bfu)^\top\big)\geq \ALPHA
I_{m\ti m}.
\end{align*}
Using this, it is shown in \cite{MieSch21?ESPS} that $\bfA$ is monotone if and
only if 
\[
(3{-}\sqrt8\,)d_3 < d_j < (3{+}\sqrt8\,) d_3 \quad \text{ for }j=1,\ 2 .
\]  
Hence, the vector-valued version of the existence theorem for similarity
profiles can be applied and for all limits $\bfU_-$ and $\bfU_+$ there
exists a unique similarity profile
$\bfU:\R \to \R^2$ connecting $\bfU_-$ and $\bfU_+$. These solutions give rise
to similarity profiles $\bfC = \Psi{\circ} \bfU$ connecting $\Psi(\bfU_-)$ and
$\Psi(\bfU_+)$ if and only if $\bfU(y) \in \mfU={[0,\infty[}^2$ for all $y \in \R$, thus
providing $\bfC(y)=\Psi(\bfU(y)) \in \mfC={[0,\infty[}^3$. In general, it seems
to be difficult to guarantee this condition, but defining $
\ol\bfu_{(\pm)}:\R\to \R^2$ via 
\[
\ol\bfu_{(\pm)}(y)=\bfU_\pm \text{ for }\pm y>0 \text{ and } 
\ol\bfu_{(\pm)}(0)=\frac12(\bfU_-{+}\bfU_+),
\]
one can show the uniform estimate 
\[
\big| \bfU(y) - \ol\bfu_{(\pm)}(y) \big| \leq C_* |\bfU_+{-} \bfU_-|,
\]
where $C_*$ only depends on $d_1$, $d_2$, and $d_3$, but not on $\bfU_\pm$. 
Thus, we obtain valid similarity profiles if $|\bfU_+{-} \bfU_-|$ is
sufficiently small compared to the distance of $\bfU_+$ and $\bfU_-$ from the
boundary of $\mfU$. In that case, similarity profiles $\bfC:\R \to \R^3$ solving
\eqref{eq:ProfEqn3Spec} exist and are unique.

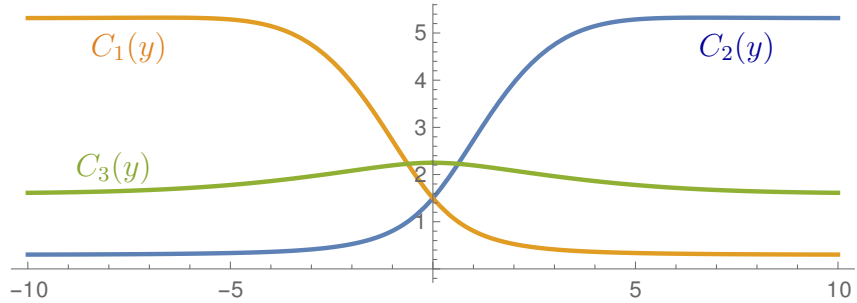
\begin{figure}[h]
\centering 
\begin{tikzpicture}
\node (PP) at (0,0){\includegraphics[width=0.7\linewidth]{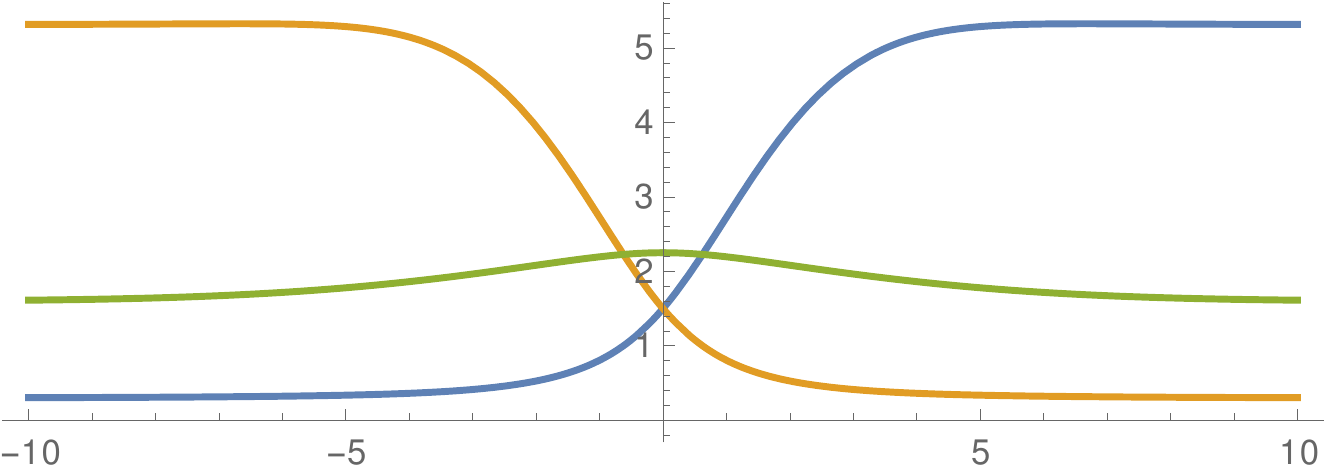}};
\node[color=brown!50!orange] at (-4,1.4) {$C_1(y)$};
\node[color=blue!60!black] at (4,1.4) {$C_2(y)$};
\node[color=green!30!brown] at (-4.2,-0.2) {$C_3(y)$};
\end{tikzpicture}

\caption{Solution $\bfC=(C_1(y),C_2(y),C_3(y))$ of \eqref{eq:ProfEqn3Spec} for
  $d_1=d_2=2$ and $d_3=10$ with limiting values $C_-\approx(5.3,0.3,1.6)$ and 
$C_+\approx(0.3,5.3,1.6) $. This symmetric solution was obtained by starting with
$\bfC(0)=(1.5,1.5,2.25)$ and $\bfC'(y)=(-1,1,0)$.}
\label{fig:Nonmonotone}
\end{figure}

In the present example we obtain nonmonotone profiles $\bfC:\R\to
\mfC\subset\R^3$. For this, consider the case $d_1=d_2$ and the limits 
\[
C_-=(A,B,AB)^\top \quad \text{and} \quad
C_+= (B,A,AB)^\top\quad \text{ with } A\neq B. 
\]
Our uniqueness result and the reflection symmetries $x \to -x$ and
$(c_1,c_2) \to (c_2,c_1)$ imply that the stationary profile $\bfC$ satisfies
$C_1(y)=C_2(-y)$ and $C_3(y)=C_3(-y)$. Using $C_1(y)C_2(y)=C_3(y)$ for all
$y\in \R$ we see that $C_3$ cannot be constant, hence it must be
nonmonotone. Figure \ref{fig:Nonmonotone} shows a corresponding example.

An interesting question is whether there is a stationary profile $\bfC$
connecting the limiting cases 
\[
C_-=\Psi(1,0)=(1,0,0)^\top  \quad \text{and} \quad
C_-=\Psi(0,1)=(0,1,0)^\top .
\]
The profile would see only one of the species $X_1$ or $X_2$ in the reservoirs
at $\pm\infty$, however in the middle region all three species must be present to allow
the generation of the other species.

\subsection{Two reactions for three species}
\label{su:TwoReact3Spec} 

Consider the two reactions $2X_1\leftrightharpoons X_2$ and  $  X_2
\leftrightharpoons X_3$ giving 
\begin{equation}
  \label{eq:2React3Spec}
  \pl_\tau \bfc = \bfD \pl_{y}^2 \bfc
  -k_1\big(c_1^2-c_2\big)\ThreeVect2{-1}0 
  - k_2  \big(c_2-c_3\big) \ThreeVect01{-1}. 
\end{equation}
The set of equilibria is the one-parameter family given by 
\[
\bigset{\bfc\in \mfC}{ \bfR(\bfc)=0} = \bigset{(A,A^2,A^2)}{ A\geq 0}.
\]
Note that the RDS system has invariant regions of the form
$\Sigma:= [b,B]\ti [b^2,B^2] \ti [b^2,B^2]$ for arbitrary
$0\leq b < B <\infty$. This means that any solution satisfying
$\bfc(0,x)\in \Sigma$ for all $x\in \R$ also satisfies $\bfc(t,x)\in \Sigma$
for all $t>0$ and $x \in \R$, see \cite{Smol94SWRD} for the theory of invariant
regions for RDS. Thus, a similarity profile connecting $C_-=(b,b^2,b^2)$ and
$C_+=(B,B^2,B^2)$ is expected to lie in the invariant region $\Sigma$.

The stoichiometric matrix is $\bfQ = (1\ \ 2 \ \ 2) \in \R^{1\ti 3}$ and 
\[
u= \bfQ \bfc = c_1 {+} 2 c_2{+} 2 c_3 \ \text{ yields } \
\Psi(u)=\ThreeVect{\sigma(u)}{(u{-}\sigma(u))/4}{(u{-}\sigma(u))/4}
\]
with $ \sigma(u)=(\sqrt{1{+}16\,u}-1)/8 $. 
With $\sigma'(u)=1/\sqrt{1{+}16\,u}\in [0,1]$ we easily see that all mappings $u
\mapsto \Psi_j(u)$ are strictly increasing such that $A(u)=
\bfQ\,\bfD\, \Psi(u)$ satisfies
\begin{align*}
&A(u)=  \tfrac{d_2{+}d_3}2 \,u + \big( d_1{-}
\tfrac{d_2{+}d_3}2\big) \,\sigma(u) \quad \text{and}
\\ 
&\min\big\{ d_1, \tfrac{d_2{+}d_3}2\big\} \leq  A'(u) \leq 
\max\big\{ d_1, \tfrac{d_2{+}d_3}2\big\}.
\end{align*}
Thus, the scalar existence theory provides for 
$0\leq U_- \leq U_+ <\infty$ a unique similarity profile $U\in
\rmC^\infty(\R; [U_-,U_+])$ that is strictly increasing.  

As a consequence, the profile equation  
\begin{equation}
  \label{eq:ProfEqn3S2R}
\begin{aligned}  
& \bfD \,\bfC'' + \frac y2  \bfC' + \Lambda_1\ThreeVect2{-1}0
  +\Lambda_2 \ThreeVect01{-1} =0, \\ 
&    C_1^2=C_2=C_3 \ \text{ and } \
    \bfC(\pm\infty) = \ThreeVect{B_\pm}{B_\pm^2}{B_\pm^2}
\end{aligned}
\end{equation}
has for all $B_-\leq B_+$ a unique solution $\bfC$ and each component $C_j$ is
strictly increasing, and hence lying in the invariant region
$\Sigma=[B_-,B_+]\ti [B_-^2,B_+^2] \ti [B_-^2,B_+^2]$. 

In this example we have the three diffusion fluxes $Q_j^\mafo{diff}(y)=-
d_jC'_j(y)$ for the three species $X_j$ and two reaction fluxes
$\Lambda^\mafo{react}_1$ and $\Lambda^\mafo{react}_2$ for the reactions 
$2X_1\leftrightharpoons X_2$ and  $  X_2
\leftrightharpoons X_3$, respectively.

\section{\label{se:GinzLand} Diffusive mixing of roll pattern} 

For a complex-valued amplitude $A(t,x) \in \C$ the real 
Ginzburg-Landau equation (i.e. the coefficients are real)
\begin{equation}
  \label{eq:RGLeqn}
  \dot A = A_{xx} + A - |A|^2 A
\end{equation}
is an important model in bifurcation theory and pattern formation. The equation
appears as amplitude or envelope equation in many partial differential
equations\cite{KiScMi92VMEE,Eckh93GLMA,Schn94EGLA,Miel02GLER,Miel15DAEE} as
well as delay equations with large delay\cite{WYHS10CDDD, YLWM15SAES}.

It has an
explicit two-parameter family of steady state pattern  in form of the  
role solutions $A(x)=U_{\eta,\varphi}(x):=\sqrt{1{-}\eta^2} \:\ee^{\ii (\eta
  x+\varphi)}$ with wave 
number $\eta \in [-1,1]$ and phase $\varphi\in [0,2\pi]$.

Starting from \cite{BriKup92RGGL,ColEck92SPSG}, it was shown in
\cite{GalMie98DMSS} that asymptotically self-similar profiles exist that
connect two different role solutions $U_{\eta_-,\varphi_-}$ at $x\to -\infty$
and $U_{\eta_+,\varphi_+} $ at $x\to \infty$. Indeed, the monotone operator
approach for showing the existence of self-similar profiles
was initiated there, see Theorem 3.1 in \cite{GalMie98DMSS} and further
developed in \cite{MieSch21?ESPS}.

Writing $A=r\ee^{\ii u}$ and assuming $r(t,x)>0$ the real Ginzburg-Landau
equation can be rewritten as the coupled system $\dot r = r_{xx}  + r\big(
1{-}r^2 {-} u_x^2\big)$ and $\dot u = u_{xx} + 2r_x u_x/r$.

Following Sec.\,2 in \cite{GalMie98DMSS} we transform the system into scaling
variables via $ t = \ee^\tau$, \ $x = \ee^{\tau/2} y$,
\begin{align*}
 \ \psi(\tau,y) = \ee^{-\tau/2} u(\ee^\tau,
 \ee^{\tau/2} y), \ \text{ and }  \ \rho(\tau,y) = r(\ee^\tau, \ee^{\tau/2} y). 
\end{align*} 
Note that $u $ and $\psi$ are related with an additional factor $\ee^{\tau/2}$,
which is necessary to match the linear behavior $u(t,x) \approx c_\pm +
\eta_\pm x$ for $x \to \pm \infty$. With this definition we still have 
$\psi(\tau, y) \approx c_\pm \ee^{-\tau/2} +
\eta_\pm y$ for $y \to \pm \infty$.

The transformed system reads 
\begin{align*}
\psi_\tau& = \psi_{yy} + \frac y2\, \psi_y - \frac12\, \psi + 2 \frac{\rho_y}
\rho\, \psi_y, 
\\
\rho_\tau&= \rho_{yy} + \frac y2\, \rho_y  + \ee^\tau
\,\rho\,(1{-}\rho^2{-} \psi_y^2\big) . 
\end{align*}
Thus, we see that for $\tau\gg 1$ we have the relation $\rho^2+ \psi_y^2
\approx 1$. Inserting the constraint $\rho=\big( 1{-}\psi_y^2\big)^{1/2}$ we
obtain the scaled phase-diffusion equation 
\begin{equation}
  \label{eq:psi.Eqn}
  \psi_\tau = \big( \Phi(\psi_y)\big)_y + \frac y2\, \psi_y - \frac12\, \psi \ 
\text{ with }\Phi'(\eta) = \frac{ 1 {-}3\eta^2 }{ 1{-}\eta^2}  .
\end{equation}
Moreover, the limiting equation for $\rho$ reads 
\[
\rho_\tau = \rho_{yy} + \frac y2\, \rho_y  +\Lambda, \quad \rho^2+\psi_y^2=1,
\]
where in principle it is possible to determine the Lagrange multiplier
$\Lambda$ from the constraint $\rho^2+\psi_y^2=1$ and \eqref{eq:psi.Eqn}. 

Using $\eta(\tau,y)$ and differentiation once we obtain a scaled diffusion
equation like the PME:
\begin{equation}
  \label{eq:eta.Eqn}
  \eta_\tau = \big( \Phi(\eta)\big)_{yy} + \frac y2\, \eta_y = 
\big( \Phi'(\eta) \eta_y\big)_y + \frac y2\, \eta_y,
\end{equation}
which shows that the equation is well-posed only for $\Phi'(\eta)>0$, i.e.\ 
$|\eta|<1/\sqrt3$, where $|\eta|>1/\sqrt3$ leads to the celebrated Eckhaus
instability \cite{Eckh65SNLS,EcGaWa95PSEI,Miel97ISRS}. 

For all $\eta_-,\eta_+\in {]{-}1/\sqrt3, 1/\sqrt3[}$ there exists a unique steady
profile $\ol \eta$ for \eqref{eq:eta.Eqn} and via
\[
\ol\psi(y) = \eta_- y + \int_{-\infty}^y \!\ol\eta(s){-}\eta_- \dd s = 
\eta_+ y - \int_{y}^\infty\! \ol\eta(s){-}\eta_+ \dd s
\]
we obtain the steady profile $\ol\psi$ for \eqref{eq:psi.Eqn} with the correct
asymptotics for $x \to \pm \infty$. 
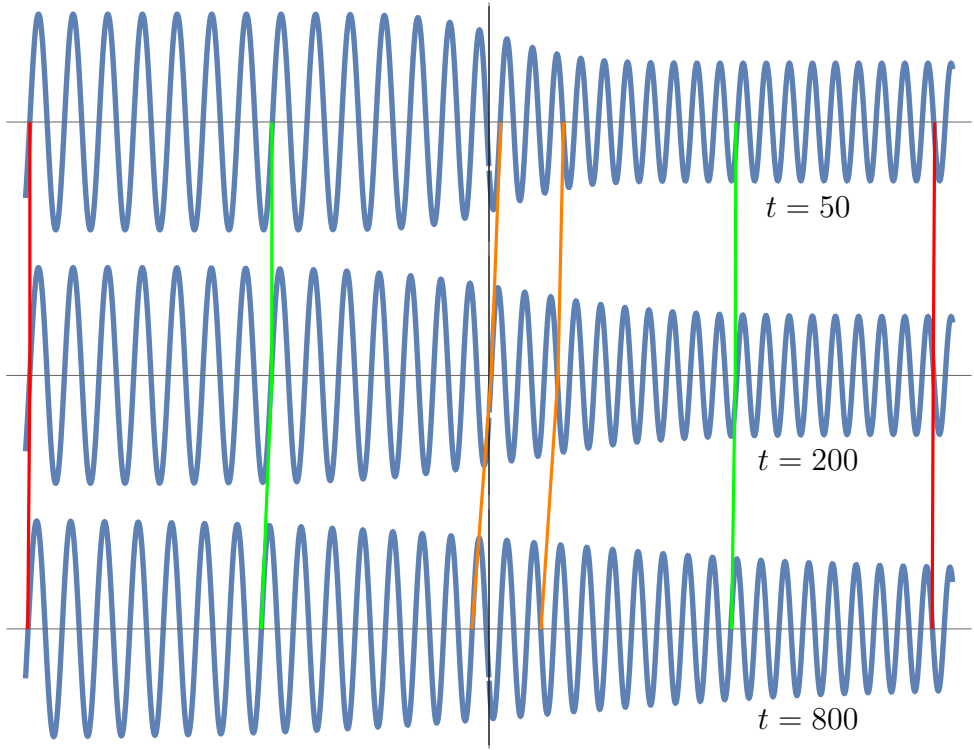
\begin{figure}
\centering 
\begin{tikzpicture}[scale=1.4]
\node at (0,4.8){\includegraphics[width=0.8\linewidth]{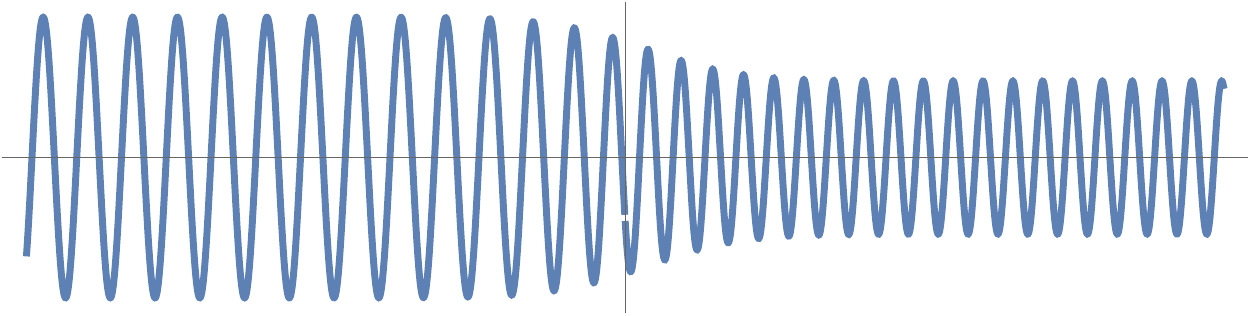}};
\node at (0,2.4){\includegraphics[width=0.8\linewidth]{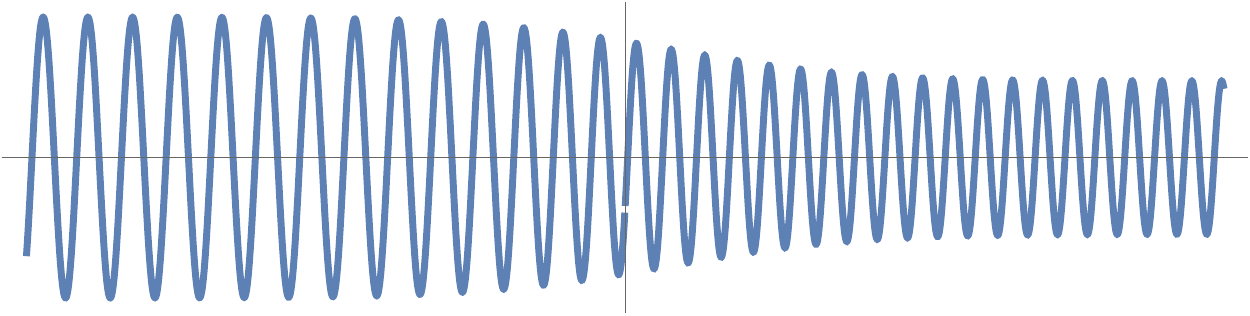}};
\node at (0,0){\includegraphics[width=0.8\linewidth]{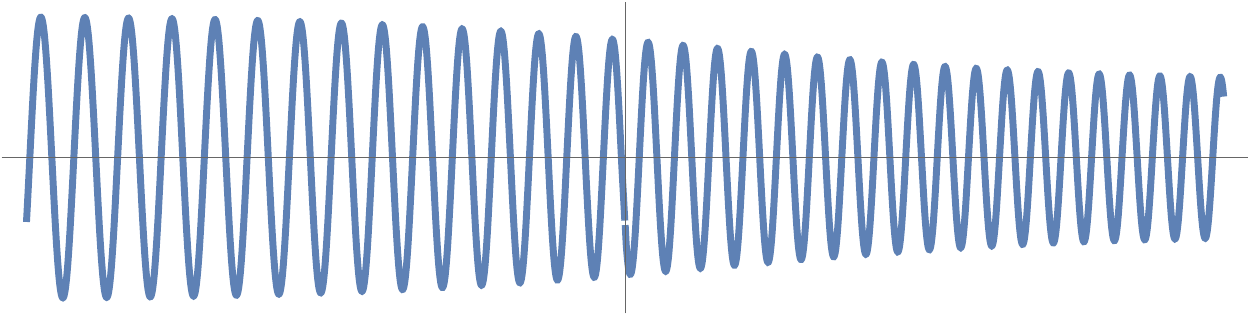}};

\node at (3,4) {$t=50$};
\node at (3,1.6) {$t=200$};
\node at (3,-0.85) {$t=800$};
\draw[thin] (0,5.9)--(0,-1.1);

\draw[very thick,color=red] (4.189,4.8) -- (4.175,2.4)--(4.165,0);
\draw[very thick,color=green] (2.32,4.8) -- (2.314,2.4)--(2.28,0);
\draw[very thick,color=orange] (0.7,4.8) -- (0.648,2.4)--(0.489,0);
\draw[very thick,color=orange] (0.11,4.8) -- (0.02,2.4)--(-0.16,0);
\draw[very thick,color=green] (-2.04,4.8) -- (-2.0415,2.4)--(-2.14,0);
\draw[very thick,color=red] (-4.315,4.8) -- (-4.315,2.4)--(-4.335,0);
\end{tikzpicture}

\caption{\label{fig:GinzLand}The three graphs display $\Re\, A(t,x)$ for
$t=50,\,200,\,800$ for the case $\eta_-=0.45$ and $\eta_+=0.3$. The vertical
connections between the graphs show the motion of the zeros. It is very slow
for large $|x|$ (red lines) and is larger for $|x|$ smaller (orange lines).} 
\end{figure}

In Figure \ref{fig:GinzLand} we sketch the self-similar behavior of the
solution $A(t,x)$ for three different times for $0< \eta_0=0.3 < \eta_+=0.45$. 
The diffusive mixing leads to a motion of the zeros of $\Re\,  A(t,x)$
to the left. The speed $v(t,x)$ of the zeros located at $x$ at time $t$ follows a self-similar profile, namely 
\[
v(t,x)= \frac1{\sqrt{1{+}t}} \,V\big( x/\sqrt{1{+}t}\big) \ \text{ with } 
V(y) = \frac y2 - \frac{\ol\psi(y)}{2 \ol\psi'(y)}.
\]
Here $V$ can be calculated by observing that a zero placed at $x_0$ for time
$t=0$ corresponds to the phase $u_0=\ol\psi(x_0)$. As the phase evolves like
$u(t,x)= \sqrt{1{+}t}\; \ol\psi\big( x/\sqrt{1{+}t} \big)$, the position of the
chosen zero has the form
$x(t)=\sqrt{1{+}t}\; H\big(\ol\psi(x_0)/\sqrt{1{+}t}\big)$, where $H$ is the
inverse mapping of $\ol\psi$. Taking the time
derivative and transforming back, provides the result.

\section{\label{sec:Concl} Conclusion}
 
In the previous sections we have shown that there are three different types of
self-similar behavior for evolution equations on $\R^d$:\medskip

(1) The
classical \emph{self-similar solutions} $\bfu(t,x)= (1{+}t)^{-\alpha} \bfU\big(
x/(1{+}t)^\beta\big)$ solve the underlying system exactly. As examples
we considered the Barenblatt solutions for the PME \eqref{eq:PME10} or the
exact solutions constructed via $\bbE$ for reaction-diffusion system in 
Section \ref{su:OneReact} in the special case $\d_1=d_2$
and $\beta=\gamma$.\medskip

(2) A slightly more general occurrence of \emph{asymptotically self-similar
  behavior} appears in Section \ref{se:CGPS} where the scaled equation is
nonautonomous with a term $\ee^{-\gamma d\tau}$ that vanishes for
$\tau\to \infty$. In such situations one can establish existence of profiles by
neglecting the term involving the decaying factor $\ee^{-\gamma d\tau}$,
determining the arising steady states (which are hopefully stable), and finally
applying a perturbation argument to obtain the convergence to the desired steady
state. This then shows that the solutions behave asymptotically self-similar.

However, we emphasize that even in the case treated in Section \ref{se:CGPS}
there is a subtle interplay between the conserved quantities. Only by the help of
the term involving $\ee^{-\gamma d\tau}$ it is possible to show that all the
initial energy $\calE(v(0),k(0))$ is finally turned into turbulent kinetic
energy.\medskip

(3) The most challenging situation occurs in the cases where the asymptotic
behavior is obtained by a constraint arising from an exponentially growing
factor $\ee^{\tau}$ that forces the system into a local equilibrium state. In
that case the natural limit problem is a \emph{constrained system} like in the RDS
case in Section \ref{se:RDS} and in the Ginzburg-Landau case in Section
\ref{se:GinzLand}. The term $\ee^\tau \bfR(\bfc)$ is of the limiting type
``$\infty\cdot \bm0$'' and needs to be replaced by a Lagrange multiplier
(possibly vector-valued, see Section \ref{su:TwoReact3Spec}).    
\medskip

In the cases (2) and (3) there remains to study the important question
whether or not the formally obtained self-similar profiles are indeed stable. 
This task is not addressed here, but first results are obtained in
\cite{VanPel77ABSN,GalMie98DMSS,Vazq07PMEM,GalSli22DREE,MieSch22?CSSP}. 
\medskip

The description of asymptotically self-similar behavior via the corresponding
similarity profiles in the scaled variables leads to a natural interpretation
of this behavior as a steady state in the sense of \emph{non-equilibrium steady
states}, because the stationarity of the system is only induced by the
renormalization of the time-dependent scaling variables. Hence, there are
nontrivial fluxes that balance the masses or energies in a suitable way. 
The major observation is that the appearing Lagrange multipliers are exactly
the missing fluxes that  are still relevant despite the fact that the system is
locally equilibrated. 

\bigskip

\paragraph*{Acknowledgment.} The research of A.M. was partially supported by
DFG via the Berlin Mathematics Research Center MATH+ (EXC-2046/1, project ID:
390685689), subproject ``DistFell''. The research of S.S. was supported by DFG
via SFB\,910 ``Control of self-organizing nonlinear systems'' (project number
163436311), subproject A5 ``Pattern formation in coupled parabolic systems''.


\newcommand{\etalchar}[1]{$^{#1}$}
\def\cprime{$'$}
\providecommand{\bysame}{\leavevmode\hbox to3em{\hrulefill}\thinspace}

\end{document}